\documentstyle[11pt]{article}
\begin{document}

\title{Transformations of quasilinear systems
originating from the projective theory of congruences}
\author{{\Large Ferapontov E.V. } \\
    Department of Mathematical Sciences \\
    Loughborough University \\
    Loughborough, Leicestershire LE11 3TU \\
    United Kingdom \\
    and \\
    Centre for Nonlinear Studies \\
    Landau Institute of Theoretical Physics\\
    Academy of Science of Russia, Kosygina 2\\
    117940 Moscow, GSP-1, Russia\\
    e-mail: {\tt fer@landau.ac.ru}
}
\date{}
\maketitle

\newtheorem{theorem}{Theorem}

\pagestyle{plain}

\maketitle

\begin{abstract}
We continue the investigation of the correspondence between systems of
conservation laws
$$
u^i_t=f^i(u)_x, ~~~~ i=1,...,n
$$
and $n$-parameter families of lines (congruences of lines) in $A^{n+1}$
defined by the equations
$$
y^i=u^i\ y^0-f^i(u).
$$
Relationship between "additional" conservation laws
$$
h(u)_t=g(u)_x
$$
and hypersurfaces conjugate to a congruence is established. This
construction allows us to introduce, in a purely geometric way, the
L\'evy transformations of semihamiltonian systems.
Correspondence between commuting flows
$$
u^i_{\tau}=q^i(u)_x, ~~~~ i=1,...,n
$$
and certain $n$-parameter families of planes
containing the lines of the congruence is pointed out. In the particular
case $n=2$ this construction provides an explicit parametrization of surfaces,
harmonic to a given congruence. Adjoint L\'evy transformations
of semihamiltonian systems are discussed.
Explicit formulae for the L\'evy and adjoint L\'evy transformations of the
characteristic velocities are set down.

A closely related construction of the Ribaucour congruences of spheres
is discussed in the Appendix.

\bigskip

Subj. Class.: Differential Geometry, Partial Differential Equations.



\end{abstract}

\section{Systems of conservation laws. Equations for the conserved
densities}

We consider systems of conservation laws
\begin{equation}
u^i_t=f^i(u)_x=v^i_j(u)  u^j_x, ~~~~~
v^i_j=\frac{\partial f^i}{\partial u^j}.
\label{cons}
\end{equation}
Eigenvalues $\lambda^i$ of the matrix $v^i_j$ are called the
characteristic velocities of system (\ref{cons}). In what follows system
(\ref{cons}) is assumed to be strictly hyperbolic, so that $\lambda^i$ are
real and pairwise distinct.
Let $\stackrel{\rightarrow}{\xi_i}=(\xi_i^1(u), ..., \xi_i^n(u))^t$ be the
corresponding right eigenvectors:
$$
v\ \stackrel{\rightarrow}{\xi_i}=\lambda^i
\stackrel{\rightarrow}{\xi_i},
~~ {\rm or, ~ in ~ the ~ components,} ~~ v^s_k\ \xi^k_i=\lambda^i\
\xi^s_i.
$$
We denote by $L_i=\xi^k_i\frac{\partial}{\partial u^k}$
the Lie derivative along the vector field
$\stackrel{\rightarrow}{\xi_i}$ and introduce commutator
expansions
$$
[L_i, L_j]=L_iL_j-L_jL_i=c^k_{ij}\ L_k,
$$
where $c^k_{ij}$ are certain functions of $u$. Let
$$
h(u)_t=g(u)_x
$$
be any conservation law of system (\ref{cons}). Its density $h$ and
flux $g$ satisfy the equations
$$
\frac{\partial g}{\partial u^k}=\frac{\partial h}{\partial u^s}\ v^s_k.
$$
Contraction with
$\stackrel{\rightarrow}{\xi_i}=(\xi_i^1, ..., \xi_i^n)^t$ results in
$$
\frac{\partial g}{\partial u^k}\ \xi^k_i=\frac{\partial h}{\partial u^s}\
v^s_k \ \xi^k_i,
$$
or
\begin{equation}
L_i g= \lambda^i L_i h, ~~~ i=1,..., n.
\label{int}
\end{equation}
Equations (\ref{int}) are the defining equations for the conserved
densities $h$ and the corresponding fluxes $g$. The compatibility
conditions of (\ref{int}) are of the form
$$
L_i(L_jg)-L_j(L_ig)=c^k_{ij}L_kg,
$$
or, taking into account (\ref{int}),
$$
L_i(\lambda^jL_jh)-L_j(\lambda^iL_ih)=c^k_{ij}\lambda^kL_kh.
$$
This results in the following linear second-order system for the conserved
densities $h$:
\begin{equation}
L_iL_jh=\frac{L_j\lambda^i}{\lambda^j-\lambda^i}\ L_ih+
\frac{L_i\lambda^j}{\lambda^i-\lambda^j}\ L_jh+
c^k_{ij} \ \frac{\lambda^i-\lambda^k}{\lambda^i-\lambda^j}\ L_kh,  ~~~ i\ne j.
\label{h}
\end{equation}
In particular, $h=u^1, ..., u^n$ satisfy (\ref{h}). It should be pointed
out that in the generic situation (to be more precise, in the case
$c^k_{ij}\ne 0$ for any $i\ne j\ne k$) the overdetermined system
(\ref{h}) possesses at most finite-dimensional linear space of solutions.
In what follows we will make use of the equations satisfied by the ratio
$\varphi =\frac{g}{h}$, which can be obtained by rewriting (\ref{int}) in
the form
$$
L_i(\varphi \ h)=\lambda^i\ L_ih,
$$
or, equivalently,
\begin{equation}
L_i \ln h=\frac{L_i\varphi}{\lambda^i-\varphi}.
\label{hphi}
\end{equation}
The compatibility conditions of (\ref{hphi}) imply the following nonlinear
second-order system for $\varphi$:
\begin{equation}
\begin{array}{c}
L_iL_j\varphi=
\left(\frac{1}{\varphi-\lambda^i}+\frac{1}{\varphi-\lambda^j}\right)
L_i\varphi \ L_j\varphi+
\frac{L_j\lambda^i}{\lambda^j-\lambda^i}
\frac{\varphi-\lambda^j}{\varphi-\lambda^i}L_i\varphi+
\frac{L_i\lambda^j}{\lambda^i-\lambda^j}
\frac{\varphi-\lambda^i}{\varphi-\lambda^j}L_j\varphi+ \\
\ \\
c^k_{ij} \ \frac{\lambda^i-\lambda^k}{\lambda^i-\lambda^j}
\frac{\varphi-\lambda^j}{\varphi-\lambda^k}
L_k\varphi.
\end{array}
\label{phi}
\end{equation}
Formula (\ref{hphi}) establishes an equivalence between the linear system
(\ref{h}) and the nonlinear system (\ref{phi}). The ratio
$\varphi=\frac{g}{h}$ naturally arises in projective differential geometry
(describing surfaces conjugate to a congruence -- see sect.2), and
in the Lie sphere geometry (parametrizing Ribaucour congruences of spheres
-- see the Appendix).

\section{Commuting flows}
System of conservation laws
\begin{equation}
u^i_{\tau}=q^i(u)_x=w^i_j(u)  u^j_x, ~~~~~
w^i_j=\frac{\partial q^i}{\partial u^j},
\label{cons1}
\end{equation}
is called the commuting flow of system (\ref{cons}), if $u_{t\tau}^i=
u^i_{\tau t}$, or, equivalently,
$$
\left(
\frac{\partial f^i}{\partial u^j} \frac{\partial q^j}{\partial u^k} u^k_x
\right)_x =
\left(
\frac{\partial q^i}{\partial u^j} \frac{\partial f^j}{\partial u^k} u^k_x
\right)_x.
$$
Equating the coefficients at $u^k_{xx}$, we arrive at the commutativity of
matrices $v=v^i_j$ and $w=w^i_j$. Thus, they possess coinciding eigenvectors
$\stackrel{\rightarrow}{\xi_i}$. Let $\mu^i$ be the characteristic velocities
of system (\ref{cons1}):
$$
w\ \stackrel{\rightarrow}{\xi_i}=\mu^i
\stackrel{\rightarrow}{\xi_i}.
$$
According to sect.1, conserved densities $h$ of system
(\ref{cons1}) satisfy the equations
\begin{equation}
L_iL_jh=\frac{L_j\mu^i}{\mu^j-\mu^i}\ L_ih+
\frac{L_i\mu^j}{\mu^i-\mu^j}\ L_jh+
c^k_{ij} \ \frac{\mu^i-\mu^k}{\mu^i-\mu^j}\ L_kh.
\label{h1}
\end{equation}
Since both systems (\ref{h}) and (\ref{h1}) share a common set of $n$
functionally independent solutions $h=u^1,..., u^n$, their coefficients must coincide
identically (if this were not the case, there would be a first-order relation between
$L_ih$, contradicting the functional independence of $u^1,..., u^n$). Thus,
\begin{equation}
\frac{L_j\mu^i}{\mu^j-\mu^i}=\frac{L_j\lambda^i}{\lambda^j-\lambda^i}
~~~~ {\rm for ~ any} ~~ i\ne j,
\label{comm1}
\end{equation}
and
\begin{equation}
c^k_{ij}\left(\frac{\mu^i-\mu^k}{\mu^i-\mu^j}-
\frac{\lambda^i-\lambda^k}{\lambda^i-\lambda^j}\right)=0
~~~~ {\rm for ~ any} ~~ i\ne j\ne k.
\label{comm2}
\end{equation}
In this form equations governing commuting flows of system (\ref{cons})
have been set down in \cite{Sevennec}.

If $n=2$, equations (\ref{comm2}) are absent. Let us consider the case
$n=3$ and assume that at least one of the coefficients $c^k_{ij}$
(with three distinct indices $i, j, k$) is nonzero. Then equations
(\ref{comm2}) imply
$$
\mu^i=\lambda^i b-a
$$
for appropriate $b$ and $a$. Substitution of this representation in
(\ref{comm1}) implies, however, that $a$ and $b$ must be constants,
so that the
commuting flow is trivial.  Hence, for $n=3$, only systems with
zero $c^k_{ij}$ (for distinct $i, j, k$) may possess nontrivial commuting
flows.

Similarly, in the case $n\geq 3$, the presence of "sufficiently many"
nonzero coefficients $c^k_{ij}$ prevents the existence of nontrivial
commuting flows.

\section{Diagonalizable systems of conservation laws}

Let us assume that all coefficients  $c^k_{ij}$ (with distinct $i, j, k$)
are  zero. In this case one can normalise eigenvectors
$\stackrel{\rightarrow}{\xi_i}$  in such a way that the Lie derivatives
$L_i$ will pairwise commute: $[L_i, L_j]=0$, so that the remaining
coefficients $c^j_{ij}$ will also be zero. The commutativity of $L_i$
implies the existense of the coordinates $R^1(u), ..., R^n(u)$,
such that $L_i$ become partial derivatives:
$$
L_i=\partial_i = \partial / \partial R^i.
$$
In the coordinates $R^i$ equations (\ref{cons}) assume  diagonal form
\begin{equation}
R_t^i=\lambda^i(R)\ R^i_x, ~~~~ i=1,...,n.
\label{riemann}
\end{equation}
Variables $R^i$ are called the Riemann invariants of system (\ref{cons}).
Systems (\ref{cons}), possessing Riemann invariants, are called
diagonalizable. Let
$$
u_t=f_x
$$
be a conservation law of system (\ref{riemann}). In the diagonalizable case
equations (\ref{int}) assume the form
$$
\partial_i f=\lambda^i\partial_i u, ~~~~ i=1,..., n,
$$
while system (\ref{h}) for the conserved densities $u$ simplifies to
\begin{equation}
\partial_i\partial_j u=a_{ij}\ \partial_i u+a_{ji}\ \partial_j u, ~~~ i\ne j,
\label{u}
\end{equation}
where $a_{ij}=\frac{\partial_j\lambda^i}{\lambda^j-\lambda^i}$.
The compatibility conditions of system (\ref{u}) are of the form
\begin{equation}
\partial_ka_{ij}=a_{ik}\ a_{kj}+a_{ij}\ a_{jk}-a_{ij}\ a_{ik}, ~~~~
i\ne j\ne k;
\label{semiham}
\end{equation}
they must be identically satisfied if we want system (\ref{u}) to possess
$n$ functionally independent solutions $u=u^1,..., u^n$. In fact, conditions
(\ref{semiham}) imply the existence of infinitely many conservation laws
parametrized by $n$ arbitrary functions of one variable. Systems (\ref{riemann})
satisfying (\ref{semiham}) are called semihamiltonian. We refer to
\cite{Tsarev}, \cite{DN}, \cite{Sevennec} for further information concerning
integrability, differential geometry and applications of semihamiltonian systems
of conservation laws. Semihamiltonian systems possess infinitely many
commuting flows
$$
R^i_{\tau}=\mu^i \ R^i_x
$$
with the characteristic velocities $\mu^i$  governed by the equations
$$
\frac{\partial_j\mu^i}{\mu^j-\mu^i}=
\frac{\partial_j\lambda^i}{\lambda^j-\lambda^i}=a_{ij},
~~~~  i\ne j.
$$

We point out that  any semihamiltonian system possesses
infinitely many different conservative representations.

\section{Systems of conservation laws and line congruences. Hypersurfaces
conjugate to a congruence}

With any system of conservation laws
$$
u^i_t=f^i(u)_x
$$
we associate an $n$-parameter family of lines
\begin{equation}
\begin{array}{c}
y^1=u^1\ y^0-f^1(u), \\
...................
\ \\
y^n=u^n\ y^0-f^n(u),
\end{array}
\label{cong}
\end{equation}
in the $(n+1)$-dimensional space $A^{n+1}$ with the coordinates
$y^0, y^1, ..., y^n$. We refer to \cite{Fer1}, \cite{Fer2} for
motivation and the most important properties of this correspondence. In
the case $n=2$ we obtain two-parameter family of lines, or a congruence of
lines in $A^3$. From the beginning of the 19th century
theory of congruences was one of the most popular chapters
of classical differential geometry -- see e.g. \cite{Finikov}.
We keep the name "congruence" for  $n$-parameter family of lines
(\ref{cong}). Any congruence possesses $n$ focal hypersurfaces
$\stackrel{\rightarrow}{{\bf r}_i}, ~~ i=1, ..., n$, ~ with the parametric
equations
$$
\stackrel{\rightarrow}{{\bf r}_i}=(y^0, y^1, ..., y^n)=
\left(\lambda^i, \ u^1 \lambda^i-f^1,
..., u^n \lambda^i-f^n\right);
$$
here $\lambda^1, ..., \lambda^n$ are the characteristic velocities of
system (\ref{cons}) -- see \cite{Fer1}, \cite{Fer2}. A line (\ref{cong})
is tangent to $\stackrel{\rightarrow}{{\bf r}_i}$ in the point with
$y^0=\lambda^i$. Let us consider a hypersurface $M^n$ with  the
radius-vector $\stackrel{\rightarrow}{{\bf r}}$ parametrized as follows:
\begin{equation}
\stackrel{\rightarrow}{{\bf r}}=(y^0, y^1, ..., y^n)=
\left(\varphi, \ u^1 \varphi-f^1,
..., u^n \varphi-f^n\right);
\label{conjugate}
\end{equation}
here $\varphi (u)$ is an arbitrary function which is assumed to be
different from $\lambda^i$ so that $M^n$ is not focal. A line (\ref{cong})
meets $M^n$ in the point with $y^0=\varphi$. Obviously, any
hypersurface $M^n\in A^{n+1}$ can be parametrized in the form
(\ref{conjugate}) for an appropriate function $\varphi$. We say that
hypersurface $M^n$ is conjugate to congruence (\ref{cong}) if and only
if
$$
L_iL_j\stackrel{\rightarrow}{{\bf r}} \in TM^n ~~~ {\rm for \ any}
~ i \ne j.
$$
Geometrically, this means that the developable surfaces
of congruence (\ref{cong}) meet $M^n$ in the curves of a
conjugate net.  In 3-space the notion of  conjugacy between a surface
and a congruence was introduced by Guichard (see \cite{Eisenhart}, chapter
1; \cite{Finikov}).

\begin{theorem}

Hypersurface (\ref{conjugate}) is conjugate to a congruence
if and only if $\varphi$
is representable in the form $\varphi=\frac{g}{h}$, where
$h_t=g_x$ is a conservation law of system (\ref{cons}).

\end{theorem}

\centerline{Proof:}

The tangent space of $M^n$ is spanned by the vectors
\begin{equation}
L_j\stackrel{\rightarrow}{\bf r}=(L_j\varphi)\stackrel{\rightarrow}{\bf
U}+(\varphi-\lambda^j)L_j\stackrel{\rightarrow}{\bf U},
\label{tangent}
\end{equation}
where $\stackrel{\rightarrow}{\bf U}$ denotes the $(n+1)$-vector $(1,
u^1, ..., u^n)$. Hence,
\begin{equation}
L_j\stackrel{\rightarrow}{\bf U}=\frac{L_j\varphi}{\lambda^j-\varphi}
\stackrel{\rightarrow}{\bf U} ~~ {\rm mod} ~ TM^n.
\label{tm}
\end{equation}
Let us compute $L_iL_j\stackrel{\rightarrow}{\bf r}$:
$$
L_iL_j\stackrel{\rightarrow}{\bf r}=(L_iL_j\varphi)
\stackrel{\rightarrow}{\bf U} +(L_j\varphi)L_i\stackrel{\rightarrow}{\bf
U}+L_i(\varphi-\lambda^j)L_j\stackrel{\rightarrow}{\bf
U}+(\varphi-\lambda^j)L_iL_j\stackrel{\rightarrow}{\bf U}.
$$
Inserting here $L_iL_j\stackrel{\rightarrow}{\bf U}$ from (\ref{h})
and keeping in mind (\ref{tm}), we arrive at
$$
\begin{array}{c}
L_iL_j\stackrel{\rightarrow}{\bf r}=
(L_iL_j\varphi+L_j\varphi
\frac{L_i\varphi}{\lambda^i-\varphi}+L_i(\varphi-\lambda^j)
\frac{L_j\varphi}{\lambda^j-\varphi}+ \\
\ \\
(\varphi-\lambda^j)(
\frac{L_j\lambda^i}{\lambda^j-\lambda^i}\frac{L_i\varphi}{\lambda^i-\varphi}
+\frac{L_i\lambda^j}{\lambda^i-\lambda^j}\frac{L_j\varphi}{\lambda^j-\varphi}
+c^k_{ij}
\frac{\lambda^i-\lambda^k}{\lambda^i-\lambda^j}
\frac{L_k\varphi}{\lambda^k-\varphi} ) )
\stackrel{\rightarrow}{\bf U} ~~{\rm mod} ~ TM^n.
\end{array}
$$
Hence, $L_iL_j\stackrel{\rightarrow}{\bf r}\in TM^n$ if and
only if the coefficient at $\stackrel{\rightarrow}{\bf U}$ vanishes.
The resulting system for $\varphi$ identically coincides with (\ref{phi}).

Thus, hypersurfaces conjugate to a congruence (\ref{cong}) are parametrized by
conservation laws of system (\ref{cons}). According to \cite{Eisenhart},
two hypersurfaces conjugate to one and the same congruence are said to be in
relation F (or related by a Fundamental transformation).

{\bf Remark 1}. The case $\varphi=\lambda^i$ requires a special
treatment. In this case $M^n$ coincides with the i-th focal hypersurface
of a congruence. A direct computation shows that the i-th focal
hypersurface is conjugate to a congruence if and only if $c^i_{jk}=0$
for any $j, k\ne i$ (i is fixed!). This is equivalent to the existence of
a function $R^i(u)$ (called the i-th Riemann invariant) such that
$$
R^i_t=\lambda^i\ R^i_x;
$$
in particular, all focal hypersurfaces are conjugate to a congruence if
and
only if system (\ref{cons}) possesses $n$ Riemann invariants. The
proof and some further details can be found in \cite{Fer2}, see also
\cite{Akivis}.

{\bf Remark 2}. If conservation law $h_t=g_x$
is a linear combination of  conservation laws (\ref{cons}),
hypersurface $M^n$ degenerates into a hyperplane
(which is automatically conjugate to any congruence). Thus, only
"additional" conservation laws give rise to nontrivial conjugate
hypersurfaces.

{\bf Remark 3}. Conjugate hypersurfaces always appear in 1-parameter
families
since, for a fixed density $h$, one can add a constant $c$
to the flux $g$. The corresponding family of conjugate hypersurfaces
$\stackrel{\rightarrow}{\bf r}_c$ determined by $\varphi _c =\frac{h}{g+c}$
forms a parallel family, that is, the directions
$L_i \stackrel{\rightarrow}{\bf r}_c$ are independent of $c$. This
immediately follows from  (\ref{tangent}), since the ratio
$\frac{L_j \varphi}{\varphi - \lambda^j}=-\frac{L_jh}{h}$ does not depend on $c$.

\section{Surfaces harmonic to a congruence}

In this section we consider 2-component systems of conservation laws
\begin{equation}
\begin{array}{c}
u_t^1=f^1_x, \\
u_t^2=f^2_x,
\end{array}
\label{cons2}
\end{equation}
and the associated congruences of lines in  $A^3$:
\begin{equation}
\begin{array}{c}
y^1=u^1\ y^0-f^1, \\
y^2=u^2\ y^0-f^2.
\end{array}
\label{cong2}
\end{equation}
Let
\begin{equation}
\begin{array}{c}
u_{\tau}^1=q^1_x, \\
u_{\tau}^2=q^2_x,
\end{array}
\label{comm3}
\end{equation}
be a commuting flow of system (\ref{cons2}). In Riemann invariants $R^1, R^2$
(we point out that any two-component system is diagonalizable) equations
(\ref{cons2}) and (\ref{comm3}) assume the forms
$$
\begin{array}{c}
R_{t}^1=\lambda^1\ R^1_x, \\
R_{t}^2=\lambda^2\ R^2_x,
\end{array}
$$
and
$$
\begin{array}{c}
R_{\tau}^1=\mu^1\ R^1_x, \\
R_{\tau}^2=\mu^2\ R^2_x,
\end{array}
$$
respectively.
The densities $u=(u^1, u^2)$ and the fluxes $f=(f^1, f^2)$,  $q=(q^1, q^2)$
satisfy the equations
$$
\partial_i f=\lambda^i\ \partial_i u, ~~~~
\partial_i q=\mu^i\ \partial_i u, ~~~~
i=1, 2.
$$
With the commuting flow (\ref{comm3}) we associate   $2$-parameter
family of planes in $A^3$ defined by the equations
\begin{equation}
\frac{y^1-u^1\ y^0+f^1}{q^1} =  \frac{y^2-u^2\ y^0+f^2}{q^2}.
\label{plane2}
\end{equation}
The family of planes (\ref{plane2}) has the following remarkable properties:

{\it
1. Each plane $\pi$ from family (\ref{plane2}) contains a line $l$ of
congruence (\ref{cong2}).

2. The congruence of lines $l_1=\pi \cap \partial_1 \pi $
is conjugate to the focal surface $\stackrel{\rightarrow}{{\bf r}_1}$
of congruence (\ref{cong2}). Similarly, the congruence of lines
$l_2=\pi \cap \partial_2 \pi $ is conjugate to
$\stackrel{\rightarrow}{{\bf r}_2}$. Lines $l_1$ and $l_2$ are called the
characteristics of plane $\pi$. Characteristic $l_1$ (resp, $l_2$),
meets the line $l$ in the point of tangency of $l$ with the focal surface
$\stackrel{\rightarrow}{{\bf r}_1}$
(resp, $\stackrel{\rightarrow}{{\bf r}_2}$). }

The proof follows from the explicit parametrization of congruences  $l_1$,
$l_2$:

\centerline{Congruence $l_1$}
$$
\begin{array}{c}
y^1=\left(u^1-\frac{q^1}{\mu^1}\right)y^0-
\left(f^1-\frac{\lambda^1q^1}{\mu^1}\right), \\
\ \\
y^2=\left(u^2-\frac{q^2}{\mu^1}\right)y^0-
\left(f^2-\frac{\lambda^1q^2}{\mu^1}\right).
\end{array}
$$

\centerline{Congruence $l_2$}
$$
\begin{array}{c}
y^1=\left(u^1-\frac{q^1}{\mu^2}\right)y^0-
\left(f^1-\frac{\lambda^2q^1}{\mu^2}\right), \\
\ \\
y^2=\left(u^2-\frac{q^2}{\mu^2}\right)y^0-
\left(f^2-\frac{\lambda^2q^2}{\mu^2}\right).
\end{array}
$$
Obviously, the line $l_1$ passes through the point
$$
(y^0, y^1, y^2)= \left(\lambda^1, \ u^1 \lambda^1-f^1,
\ u^2 \lambda^1-f^2\right)
$$
of the focal surface $\stackrel{\rightarrow}{{\bf r}_1}$.
Similarly, the line $l_2$ passes through the point
$$
(y^0, y^1, y^2)= \left(\lambda^2, \ u^1 \lambda^2-f^1,
\ u^2 \lambda^2-f^2\right)
$$
of the focal surface $\stackrel{\rightarrow}{{\bf r}_2}$.
The point of intersection $l_1\cap l_2\in \pi$ has the coordinates
\begin{equation}
\begin{array}{c}
y^0=\frac{\lambda^2 \mu^1-\lambda^1 \mu^2}{\mu^1-\mu^2}, \\
\ \\
y^1=\frac{\lambda^2 \mu^1-\lambda^1 \mu^2}{\mu^1-\mu^2}\ u^1+
\frac{\lambda^1-\lambda^2 }{\mu^1-\mu^2}\ q^1-f^1, \\
\ \\
y^2=\frac{\lambda^2 \mu^1-\lambda^1 \mu^2}{\mu^1-\mu^2}\ u^2+
\frac{\lambda^1-\lambda^2 }{\mu^1-\mu^2}\ q^2-f^2,
\end{array}
\label{harm}
\end{equation}
and sweeps a surface in $A^3$. By a
construction, surface (\ref{harm}) is the envelope of the family of planes
(\ref{plane2}). It has the following geometric properties:

{\it
1. Each tangent plane $\pi$ of  surface (\ref{harm}) contains a line $l$
of  congruence (\ref{cong2}). By a construction, $\pi$ and $l\in \pi$
correspond to the same values of  parameters $R^1, R^2$. Thus, one can speak of
the correspondence between lines (\ref{cong2}) and points of  surface
(\ref{harm}).

2. The net $R^1, R^2$ on  surface (\ref{harm}) is conjugate. In other words,
developable surfaces of   congruence (\ref{cong2}) correspond to a conjugate net
on  surface (\ref{harm}). }

Surfaces, satisfying the properties 1, 2, are called
harmonic to congruence (\ref{cong2}) -- see \cite{Finikov}, p.251.
Formulae (\ref{harm}) provide an explicit parametrization of surfaces, harmonic
to congruence (\ref{cong2}), by commuting flows of system (\ref{cons2}).
Conversely, any surface harmonic to congruence (\ref{cong2}) is
representable in the form (\ref{harm}).

\section{L\'evy transformations of semihamiltonian systems}

Let us consider  semihamiltonian system (\ref{riemann})
in Riemann invariants:
$$
R_t^i=\lambda^i(R)\ R^i_x, ~~~~ i=1,...,n.
$$
Conservation laws
$$
u_t=f_x
$$
of system (\ref{riemann})
satisfy the equations
$$
\begin{array}{c}
\partial_i f=\lambda^i\ \partial_i u, ~~~~~ i=1,...,n, \\
\ \\
\partial_i\partial_j u=a_{ij}\ \partial_i u+a_{ji}\ \partial_j u, ~~~ i\ne j.
\end{array}
$$
Let us choose  particular conservation law
$$
h_t=g_x
$$
of system (\ref{riemann}) and introduce new variable $U$ by the
formula
\begin{equation}
U=u-\frac{h}{\partial_{\alpha} h}\ \partial_{\alpha} u,
\label{levi}
\end{equation}
where ${\alpha}$ is fixed. Transformations of this type originate from projective
differential geometry of conjugate nets and are known as the
transformations of L\'evy \cite{Levy}, \cite{Forsyth}, p.94,
\cite{Eisenhart}, chapter 1.  In paper \cite{Doliva} transformations of
L\'evy have been identified with the vertex operators of the multicomponent KP
hierarchy. Their geometric interpretation
will be clarified in the second half of this section. We will refer to
(\ref{levi}) as to the transformation of L\'evy ${\cal L}_{\alpha}$.
A direct calculation
shows that $U={\cal L}_{\alpha}(u)$ satisfies the equations of
the same form as $u$:
\begin{equation}
\partial_i\partial_j U=A_{ij}\partial_i U+A_{ji}\partial_j U
\label{U}
\end{equation}
where the new coefficients $A={\cal L}_{\alpha}(a)$ are given by the formulae
$$
\begin{array}{c}
A_{{\alpha}i}=\left(1-\frac{a_{i{\alpha}}h}{\partial_{\alpha}h}\right)
\frac{\partial_ih}{h}, ~~~ i\ne {\alpha}, \\
\ \\
A_{ij}=a_{ij}+\partial_j
\ln \left(1-\frac{a_{i{\alpha}}h}{\partial_{\alpha}h}\right),
~~~ i\ne {\alpha}, ~~ j ~ {\rm is \ arbitrary}.\\
\end{array}
$$
Transformations ${\cal L}_{\alpha}$ can be pulled back to the transformations of the
corresponding hydrodynamic type systems:
let us introduce the system
\begin{equation}
R_T^i=\Lambda^i(R)\ R^i_X, ~~~~ i=1,...,n
\label{newriemann}
\end{equation}
with the characteristic velocities
\begin{equation}
\begin{array}{c}
\Lambda^{\alpha}=\frac{g}{h}, \\
\ \\
\Lambda^i=\frac{\displaystyle \lambda^i\partial_{\alpha}h-a_{i{\alpha}}g}
{\displaystyle \partial_{\alpha}h-a_{i{\alpha}}h}, ~~~
i\ne {\alpha}.
\end{array}
\label{newlambda}
\end{equation}

\begin{theorem}
Conservation laws
$$
U_T=F_X
$$
of system (\ref{newriemann}), (\ref{newlambda})
are the ${\cal L}_{\alpha}$-transforms of  conservation
laws
$$
u_t=f_x
$$
of system (\ref{riemann}):
$$
\begin{array}{c}
U={\cal L}_{\alpha}(u)=u-\frac{h}{\partial_{\alpha}h}\ \partial_{\alpha} u, \\
\ \\
F={\cal L}_{\alpha}(f)=f-\frac{g}{\partial_{\alpha}g}\ \partial_{\alpha} f.
\end{array}
$$
\end{theorem}

Formally, the proof of this theorem follows from the identities
$$
\partial_iF=\Lambda^i\ \partial_iU, ~~~~
A_{ij}=\frac{\partial_j \Lambda^i}{\Lambda^j-\Lambda^i},
$$
which can be verified by a direct calculation. Geometric constructions
underlying these formulae will be discussed below. System
(\ref{newriemann}), (\ref{newlambda}) will be called the 
${\cal L}_{\alpha}$-transform
of system (\ref{riemann}). Obviously, transformations  ${\cal L}_{\alpha}$
preserve the semihamiltonian property.

We also include L\'evy transformations of the Lame coefficients $h_i$
defined by the formulae
$$
\partial_j \ln h_i= a_{ij},
~~~ j\ne i.
$$
The ${\cal L}_{\alpha}$-transformed Lame coefficients are given by
$$
\begin{array}{c}
H_{\alpha}=h_{\alpha} \ \frac{h}{\partial_{\alpha}h}, \\
\ \\
H_i=h_i\ \left(1-\frac{a_{i{\alpha}}h}{\partial_{\alpha}h} \right), ~~~
i\ne {\alpha}.
\end{array}
$$
One can check directly that
$$
\partial_j \ln H_i= A_{ij},
~~~ j\ne i.
$$
L\'evy transformations of hydrodynamic type systems in Riemann invariants
are closely related to the transformations of Laplace discussed recently
in \cite{Fer3}, \cite{Kamran}.
We recall that Laplace transformation $S_{{\alpha}{\beta}}$
of system (\ref{u}) is
defined by the formula
$$
U=S_{{\alpha}{\beta}}(u)=u-\frac{\partial_{\alpha}u}{a_{{\beta}{\alpha}}},
$$
where both indices ${\alpha}\ne {\beta}$ are fixed.
Laplace transformations also induce
transformations of the characteristic velocities $\lambda^i$, the explicit
form of which has been set down in \cite{Fer3}. One can check directly
that the  L\'evy transformation ${\cal L}_{\alpha}$ of system (\ref{riemann})
is related to its L\'evy transformation ${\cal L}_{\beta}$
via the Laplace transformation $S_{{\alpha}{\beta}}$:
$$
{\cal L}_{\alpha}=S_{{\alpha}{\beta}}\circ {\cal L}_{\beta}.
$$

To clarify geometric picture underlying transformations ${\cal L}_{\alpha}$
we choose an arbitrary conservative representation
$$
u_t^i=f^i_x
$$
of system (\ref{riemann}) and introduce the associated congruence
$$
\begin{array}{c}
y^1=u^1y^0-f^1, \\
............ \\
y^n=u^ny^0-f^n.
\end{array}
$$
Let $M^n$ be hypersurface conjugate to this congruence. Following
sect.3,
we represent the radius-vector $\stackrel{\rightarrow }{\bf r}$ of $M^n$
in
the form
$$
\stackrel{\rightarrow}{\bf r}=\left(\varphi, ~ u^1\varphi - f^1, ...,
\ u^n\varphi - f^n
\right), ~~~~ \varphi = \frac{g}{h},
$$
where $h_t=g_x$ is a conservation law of system
(\ref{riemann}). Coordinate system $R^1, ..., R^n$ on
$M^n$ is conjugate, so that
$$
\partial_i\partial_j \stackrel{\rightarrow}{\bf r}\in TM^n ~~ {\rm for
\ any} ~~ i\ne j.
$$
Let
us introduce a new congruence, formed by the tangents to the
$R^{\alpha}$-curves on  hypersurface $M^n$. Parametrically, its lines
can be represented in the form
$$
\stackrel{\rightarrow}{\bf r}+t\ \partial_{\alpha}\stackrel{\rightarrow}{\bf r},
$$
or, in the components,
$$
\begin{array}{c}
y^0=\varphi+t\ \partial_{\alpha}\varphi, \\
\ \\
y^1=u^1\varphi-f^1+t\ (u^1\partial_{\alpha}\varphi+(\varphi-\lambda^{\alpha})
\partial_{\alpha}u^1), \\
............................................... \\
y^n=u^n\varphi-f^n+t\ (u^n\partial_{\alpha}\varphi+(\varphi-\lambda^{\alpha})
\partial_{\alpha}u^n).
\end{array}
$$
Inserting $t=\frac{y^0-\varphi}{\partial_{\alpha} \varphi}$
in the last n equations, we arrive at the new
congruence
\begin{equation}
\begin{array}{c}
y^1=U^1y^0-F^1, \\
................ \\
y^n=U^ny^0-F^n,
\end{array}
\label{newcongr}
\end{equation}
where
$$
\begin{array}{c}
U^1=u^1+\frac{\varphi -\lambda^{\alpha}}{\partial_{\alpha}\varphi}
\partial _{\alpha}u^1, ~~~
F^1=f^1+\varphi \frac{\varphi -\lambda^{\alpha}}{\partial_{\alpha}\varphi}
\partial_{\alpha}u^1, \\
...................................... \\
U^n=u^n+\frac{\varphi -\lambda^{\alpha}}{\partial_{\alpha}\varphi}
\partial _{\alpha}u^n, ~~~
F^n=f^n+\varphi \frac{\varphi -\lambda^{\alpha}}{\partial_{\alpha}\varphi}
\partial_{\alpha}u^n.
\end{array}
$$
Since $\frac{\varphi-\lambda^{\alpha}}{\partial_{\alpha}\varphi}
=-\frac{h}{\partial_{\alpha}h}$,
these formulae can
be rewritten in the form
$$
U=u-\frac{h}{\partial_{\alpha}h}\ \partial_{\alpha}u, ~~~
F=f-\frac{g}{\partial_{\alpha}g}\ \partial_{\alpha}f.
$$
Congruence (\ref{newcongr}) will be called the ${\cal L}_{\alpha}$-transform of the
initial congruence. The corresponding system of conservation laws
$$
U^i_T=F^i_X
$$
has the same Riemann invariants $R^1, ..., R^n$:
$$
R^i_T=\Lambda^i\ R^i_X,
$$
where $\Lambda^i$ can be computed as follows:
$\Lambda^i=\partial_iF/\partial_iU$. A direct calculation results in
formulae (\ref{newlambda}). Note that the final expressions for
$\Lambda^i$ do not depend on the particular conservative representation
$u^i_t=f^i_x$ of system (\ref{riemann}). If, for $M^n$, we choose any of
the
focal hypersurfaces of the congruence (which are all conjugate to a
congruence if the system possesses Riemann invariants),
the above construction gives transformations of Laplace.

Formula (\ref{U}) shows that the density $u=h$ belongs to
the kernel of the L\'evy
transformation ${\cal L}_{\alpha}$. Nevertheless, transformations
${\cal L}_{\alpha}$ can be explicitely
inverted, as we will demonstrate  in the next section.

Let us conclude with the formula for the composition of the L\'evy
transformations
$$
{\cal L}={\cal L}_n\circ ... \circ {\cal L}_2\circ {\cal L}_1
$$
corresponding to $n$ particular linearly independent conservation
laws $h^i_t=g^i_x, ~~ i=1,...,n$ of system (\ref{riemann}).
The composition is understood as follows. Let $u_t=f_x$ be an arbitrary
conservation law of system (\ref{riemann}). First of all,
we apply to  $u_t=f_x$  transformation ${\cal L}_1$, corresponding to the first
conservation law  $h^1_t=g^1_x$. Secondly, we apply to the result of the
first step  transformation ${\cal L}_2$, corresponding to the ${\cal L}_1$-transform
of  conservation law $h^2_t=g^2_x$. Proceeding in this way, we obtain
the ${\cal L}$-transformed density $U={\cal L}(u)$ and the flux $F={\cal L}(f)$ in the
following compact form:
\begin{equation}
U=\frac{det \left(
\begin{array}{cccc}
u & \partial_1u & ... & \partial_nu \\
h^1 & \partial_1 h^1 & ... & \partial_n h^1 \\
\  \\
h^n & \partial_1 h^n & ... & \partial_n h^n
\end{array}
\right)}
{det \left(
\begin{array}{ccc}
\partial_1 h^1 & ... & \partial_n h^1 \\
\  \\
\partial_1 h^n & ... & \partial_n h^n
\end{array}
\right)}, ~~~~~
F=\frac{det \left(
\begin{array}{cccc}
f & \partial_1f & ... & \partial_nf \\
g^1 & \partial_1 g^1 & ... & \partial_n g^1 \\
\  \\
g^n & \partial_1 g^n & ... & \partial_n g^n
\end{array}
\right)}
{det \left(
\begin{array}{ccc}
\partial_1 g^1 & ... & \partial_n g^1 \\
\  \\
\partial_1 g^n & ... & \partial_n g^n
\end{array}
\right)}.
\label{UF}
\end{equation}
Geometrically, the composition  ${\cal L}_n\circ ... \circ {\cal L}_2\circ {\cal L}_1$
corresponds to the following construction (compare with \cite{Finikov},
p.255-256): choose an arbitrary conservative representation
$$
u^i_t=f^i_x
$$
of system (\ref{riemann}) and introduce the corresponding congruence
(\ref{cong}):
$$
y^i=u^iy^0-f^i.
$$
Let $M_i, ~i=1, ...,n$, be $n$ hypersurfaces conjugate to congruence
(\ref{cong}). According to sect.2, they are parametrized by
$n$ particular conservation laws $h^i_t=g^i_x$ of system (\ref{riemann}).
Let $TM_i$ be the tangent hyperplanes of hypersurfaces $M_i$ in the
points of intersection with  line (\ref{cong}).
The intersection
$$
TM_1\cap ... \cap TM_n
$$
defines a new line
$$
y^i=U^iy^0-F^i;
$$
one can check directly, that the formulae for $U=U^i$ and $F=F^i$
coincide with (\ref{UF}).

\section{The adjoint transformations of L\'evy}

We again  consider  semihamiltonian systems (\ref{riemann})
$$
R_t^i=\lambda^i(R)\ R^i_x
$$
with conservation laws
$$
u_t=f_x
$$
satisfying the equations
$$
\begin{array}{c}
\partial_i f=\lambda^i\ \partial_i u, \\
\ \\
\partial_i\partial_j u=a_{ij}\ \partial_i u+a_{ji}\ \partial_j u,
\end{array}
$$
where $a_{ij}=\frac{\partial_j\lambda^i}{\lambda^j-\lambda^i}$.
Let
\begin{equation}
R^i_{\tau}=\mu ^i(R)\ R^i_x
\label{mu}
\end{equation}
be a commuting flow of system (\ref{riemann}):
$$
\frac{\partial_j \mu^i}{\mu^j-\mu^i}=a_{ij}.
$$
Let $q$ be the flux of density $u$, corresponding to this commuting flow:
$$
u_{\tau}=q_x.
$$
The flux $q$ and the density $u$ satisfy the equations
$$
\partial_i q=\mu^i\ \partial_i u.
$$
Let us introduce  new variable $U$ by the
formula
\begin{equation}
U=u-\frac{q}{\mu^{\alpha}},
\label{U1}
\end{equation}
where ${\alpha}$ is fixed. We will refer to (\ref{U1}) as to the
adjoint transformation  of L\'evy
${\cal L}^{*}_{\alpha}$.  A direct calculation shows that
$U={\cal L}^{*}_{\alpha}(u)$ satisfies the equations of the same form as $u$:
$$
\partial_i\partial_j U=A_{ij}\partial_i U+A_{ji}\partial_j U
$$
where the new coefficients $A={\cal L}^{*}_{\alpha}(a)$ are given by the formulae
$$
\begin{array}{c}
A_{{\alpha}i}=a_{{\alpha}i}+
\partial_i\ln\frac{\partial_{\alpha}\mu^{\alpha}}{\mu^{\alpha}},
~~~ i\ne {\alpha}, \\
\ \\
A_{ij}=a_{ij}+\partial_j
\ln \left(1-\frac{\mu^i}{\mu^{\alpha}}\right),
~~~ i\ne {\alpha}, ~~ j ~ {\rm is \ arbitrary}.\\
\end{array}
$$
Transformations ${\cal L}^{*}_{\alpha}$ can be pulled back to the transformations of the
corresponding hydrodynamic type systems:
let us introduce the system
\begin{equation}
R_T^i=\Lambda^i(R)\ R^i_X, ~~~~ i=1,...,n
\label{newriemann1}
\end{equation}
with the characteristic velocities
\begin{equation}
\begin{array}{c}
\Lambda^{\alpha}=
\frac{\displaystyle \lambda^{\alpha}\partial_{\alpha}\mu^{\alpha}-
\mu^{\alpha}\partial_{\alpha}\lambda^{\alpha}}
{\displaystyle \partial_{\alpha}\mu^{\alpha}}, \\
\ \\
\Lambda^i=\frac{\displaystyle \lambda^i\mu^{\alpha}-\lambda^{\alpha}\mu^i}
{\displaystyle \mu^{\alpha}-\mu^i}, ~~~
i\ne {\alpha}.
\end{array}
\label{newlambda1}
\end{equation}

\begin{theorem}
Conservation laws
$$
U_T=F_X
$$
of system (\ref{newriemann1}), (\ref{newlambda1})
are the ${\cal L}^{*}_{\alpha}$-transforms of  conservation
laws
$$
u_t=f_x
$$
of system (\ref{riemann}):
$$
\begin{array}{c}
U={\cal L}^{*}_{\alpha}(u)=u-\frac{q}{\mu^{\alpha}}, \\
\ \\
F={\cal L}^{*}_{\alpha}(f)=f-\frac{\lambda^{\alpha}q}{\mu^{\alpha}}.
\end{array}
$$
\end{theorem}

Formally, the proof of this theorem follows from the identities
$$
\partial_iF=\Lambda^i\ \partial_iU, ~~~~
A_{ij}=\frac{\partial_j \Lambda^i}{\Lambda^j-\Lambda^i},
$$
which can be verified by a direct calculation. Geometric constructions
underlying these formulae will be discussed below. System
(\ref{newriemann1}), (\ref{newlambda1}) will be called the 
${\cal L}^{*}_{\alpha}$-transform
of system (\ref{riemann}). Obviously, transformations  
${\cal L}^{*}_{\alpha}$
preserve the semihamiltonian property.

We also include ${\cal L}^{*}_{\alpha}$-transforms of the Lame coefficients $h_i$
defined by the formulae
$$
\partial_j \ln h_i= a_{ij}, ~~~ j\ne i.
$$
The ${\cal L}^{*}_{\alpha}$-transformed Lame coefficients are given by
$$
\begin{array}{c}
H_{\alpha}=h_{\alpha} \ \frac{\partial_{\alpha}\mu^{\alpha}}{\mu^{\alpha}}, \\
\ \\
H_i=h_i\ \left(1-\frac{\mu^i}{\mu^{\alpha}} \right), ~~~ i\ne {\alpha}.
\end{array}
$$
One can check directly that
$$
\partial_j \ln H_i= A_{ij}, ~~~ j\ne i.
$$
Transformations ${\cal L}^{*}_{\alpha}$ and the Laplace transformations
$S_{{\alpha}{\beta}}$ satisfy
the identities
$$
{\cal L}^{*}_{\alpha}={\cal L}^{*}_{\beta}\circ S_{{\beta}{\alpha}}.
$$
To clarify  geometric picture underlying transformations ${\cal L}^{*}_{\alpha}$
we choose an arbitrary conservative representation
$$
u_t^i=f^i_x
$$
of system (\ref{riemann}) and introduce the associated congruence
$$
\begin{array}{c}
y^1=u^1y^0-f^1, \\
............ \\
y^n=u^ny^0-f^n.
\end{array}
$$
Let
$$
u^i_{\tau}=q^i_x
$$
be a commuting flow of system (\ref{riemann}) with the characteristic velocities
$\mu$, so that
$$
\partial_iq=\mu^i\ \partial_iu,
$$
(the last identity holding for any $q=q^k, ~ u=u^k$).
Let us introduce  $n$-parameter
family of 2-planes in $A^{n+1}$ defined by the equations
\begin{equation}
\frac{y^1-u^1y^0+f^1}{q^1} = ... = \frac{y^n-u^ny^0+f^n}{q^n}.
\label{plane}
\end{equation}
The family of planes (\ref{plane}) possesses the following three important properties:

{\it
1. Each plane $\pi$ from family (\ref{plane}) contains a line $l$ of the initial
congruence.

2. Each plane $\pi$ intersects the plane $\partial_i\pi$ along a line $l_i$:
$$
l_i=\pi \cap \partial_i\pi,
$$
(we point out that two planes in $A^{n+1}$ do not necessarily intersect along a
line unless $n=2$). Geometrically, this property implies that each 1-parameter
subfamily of (\ref{plane}), specified by fixing the values of $R^k, k\ne i$,
envelopes a developable surface in $A^{n+1}$.
Lines $l_i, ~ i=1,...,n$, are called the
characteristics of plane $\pi$.

3. Congruence $l_i$ is conjugate to the i-th focal hypersurface
$$
\stackrel{\rightarrow}{{\bf r}_i}=
\left(\lambda^i, \ u^1 \lambda^i-f^1,
..., u^n \lambda^i-f^n\right)
$$
of the initial congruence $l$.}

Conversely, one can show that any $n$-parameter family of 2-planes
satisfying the properties 1 -- 3 is necessarily of the form (\ref{plane})
for an appropriate commuting flow $u^i_{\tau}=q^i_x$.

Congruence $l_{\alpha}$ will be called the 
${\cal L}^{*}_{\alpha}$-transform of the initial
congruence $l$. A direct calculation shows that $l_{\alpha}$ is
representable in the form
$$
\begin{array}{c}
y^1=U^1y^0-F^1, \\
................ \\
y^n=U^ny^0-F^n,
\end{array}
$$
where
$$
\begin{array}{c}
U^1=u^1-\frac{q^1}{\mu^{\alpha}}, ~~~
F^1=f^1-\frac{\lambda^{\alpha}q^1}{\mu^{\alpha}}, \\
...................................... \\
U^n=u^n-\frac{q^n}{\mu^{\alpha}}, ~~~
F^n=f^n-\frac{\lambda^{\alpha}q^n}{\mu^{\alpha}},
\end{array}
$$
(compare with Theorem 4). Line $l_{\alpha}$ meets the focal
hypersurface  $\stackrel{\rightarrow}{{\bf r}_{\alpha}}$
in the point
$$
\left(\lambda^{\alpha}, \ u^1 \lambda^{\alpha}-f^1,
..., u^n \lambda^{\alpha}-f^n\right).
$$
The corresponding system of conservation laws
$$
U^i_T=F^i_X
$$
has the same Riemann invariants $R^1, ..., R^n$:
$$
R^i_T=\Lambda^i\ R^i_X,
$$
(in fact, this is the analytic manifestation of the above property 3),
where the transformed characteristic velocities
$\Lambda^i=\partial_iF/\partial_iU$ coincide with (\ref{newlambda1}).
Note that the final expressions for
$\Lambda^i$ do not depend on the particular conservative representation
$u^i_t=f^i_x$ of system (\ref{riemann}).

Obviously, the inverse transformation $l_{\alpha} \to l$ is the
transformation ${\cal L}_{\alpha}$
of L\'evy. Indeed, $l_{\alpha}$ is conjugate to hypersurface
$\stackrel{\rightarrow}{{\bf r}_{\alpha}}$,
while the initial congruence $l$ consists
of the $R^{\alpha}$-tangents to hypersurface
$\stackrel{\rightarrow}{{\bf r}_{\alpha}}$.
Thus, transformations of L\'evy ${\cal L}_{\alpha}$
are the inverses of ${\cal L}^{*}_{\alpha}$.
This can be demonstrated analytically as well:

Let us consider a system
$$
R^i_t=\lambda^i\ R^i_x
$$
along with its L\'evy transform ${\cal L}_{\alpha}$ defined by formulae
(\ref{newriemann}), (\ref{newlambda}). The transformed system
(\ref{newriemann}), (\ref{newlambda}) possesses commuting flow
$$
\begin{array}{c}
\mu^{\alpha}=\frac{1}{h}, \\
\ \\
\mu^i=\frac{\displaystyle a_{i{\alpha}}}
{\displaystyle a_{i{\alpha}}h-\partial_{\alpha}h}, ~~~
i\ne {\alpha},
\end{array}
$$
(which can be obtained by a shift $g\to g+1$ in formulae
(\ref{newlambda})). Applying to the transformed system
(\ref{newriemann}), (\ref{newlambda}) transformation ${\cal L}^{*}_{\alpha}$
(generated by
the above commuting flow), we return to the initial system
$$
R^i_t=\lambda^i\ R^i_x.
$$
Conversely, let us consider transformation ${\cal L}^{*}_{\alpha}$. 
The transformed system
(\ref{newriemann1}), (\ref{newlambda1}) possesses  conservation law
$$
h_T=g_X, ~~~~ h=\frac{1}{\mu^{\alpha}}, ~~
g=\frac{\lambda^{\alpha}}{\mu^{\alpha}}
$$
(which can be obtained by a shift $q\to q-1$ in formula (\ref{U1})).
Applying to (\ref{newriemann1}), (\ref{newlambda1})
transformation ${\cal L}_{\alpha}$ (generated by this particular $h$), we also return
back to the initial system.

\section{Appendix: Ribaucour congruences of spheres}

Let $M^n$ be a hypersurface in the Euclidean space $E^{n+1}$ parametrized
by coordinates $u^1, ..., u^n$.
Let $\stackrel{\rightarrow}{\bf r}$ and $\stackrel{\rightarrow}{\bf n}$ be
the radius-vector and the unit normal of $M^n$, respectively. The
Weingarten formulae
$$
\frac{\partial \stackrel{\rightarrow}{\bf n}}{\partial u^j}= w^i_j(u)
\ \frac{\partial \stackrel{\rightarrow}{\bf r}}{\partial u^i}
$$
define the so-called Weingarten (shape) operator of hypersurface
$M^n$. Its eigenvalues and eigenvectors are called the principal
curvatures and the principal directions of $M^n$, respectively. Let us
consider a hypersphere $S$ of radius $R$ and the centre
$\stackrel{\rightarrow}{\bf r}-R\stackrel{\rightarrow}{\bf n}$, which is
tangent to $M^n$ at the point $\stackrel{\rightarrow}{\bf r}$. Specifying
$R$ as
a function of $u$, we obtain $n$-parameter family of
hyperspheres (or a congruence of hyperspheres) enveloped by
hypersurface $M^n$. Let $\tilde M^n$ be the second sheet of the envelope.
Clearly, there exists a point correspondence between both
sheets $M^n$ and $\tilde M^n$: a point $p\in M^n$ corresponds to
$\tilde p \in \tilde M^n$ if $p$ and $\tilde p$ are the two points of
tangency of one and the same hypersphere from the family $S(u)$.

{\bf Definition}.
Family of hyperspheres $S(u)$ is called the family of Ribaucour if the
principal distributions of $M^n$ correspond to the principal distributions
of $\tilde M^n$.

Let us introduce the system of hydrodynamic type
\begin{equation}
u^i_t=w^i_j(u) \ u^j_x,
\label{w}
\end{equation}
where $w^i_j$ is the Weingarten operator of $M^n$. We refer to \cite{Fer4}
for the general discussion of the correspondence between hypersurfaces and
systems of hydrodynamic type. Let
$$
h(u)_t=g(u)_x
$$
be a conservation law of system (\ref{w}).

\begin{theorem}

Congruence $S(u)$ is the congruence of Ribaucour if and only if
$R(u)$ is representable in the form
$$
R(u)=\frac{h(u)}{g(u)}
$$
for some conservation law of system (\ref{w}).
\end{theorem}

In the case $n=2$ this result (stated in a somewhat different form) can be
found in \cite{Eisenhart}. It should be emphasized
that this theorem equally applies
to hypersurfaces which do not possess a curvatute-line
parametrization (for $n=2$ such parametrization is always possible). We
hope to present the details elsewhere.

\section{Acknowledgements}
I would like to thank M.V.~Pavlov for usefull discussions. This research
was partially supported by INTAS 96-0770 and the RFFI grant 99-01-00010.


\begin{thebibliography}{99}
\addcontentsline{toc}{section}{References}


\bibitem{Fer1}S.I.~Agafonov and E.V.~Ferapontov, Systems of conservation
laws from the point of view of the projective theory of congruences,
Izv. RAN, ser. mat. 60 (1996) N.6, 3-30.


\bibitem{Fer2} S.I.~Agafonov and E.V.~Ferapontov, Theory of congruences
and systems of conservation laws, to appear in J. of Math. Sci.,
translation from Itogi Nauki, VINITI, Problems of Geometry, 1999.



\bibitem{Akivis} M.A.~Akivis and V.V.~Goldberg,
Projective differential geometry of submanifolds, Math. Library, V.49,
North-Holland, 1993.

\bibitem{Doliva} A.~Doliva, M.~Ma$\tilde n$as, L.M.~Alonso, E.~Medina,
P.M.~Santini, Charged free fermions, vertex operators and classical
theory of conjugate nets, solv-int/9803015, (1998).



\bibitem{DN} B.A.~Dubrovin and S.P.~Novikov, Hydrodynamics of weakly deformed
soliton lattices. Differential geometry and Hamiltonian theory, Uspekhi Mat. Nauk
44 (1989) N6, 29-98.

\bibitem{Eisenhart} L.P.~Eisenhart, Transformations of surfaces, second ed.,
1962.


\bibitem{Fer3} E.V.~Ferapontov, Laplace transformations of hydrodynamic
type systems in Riemann invariants: periodic sequences, J. Phys. A
30 (1997) N.19, 6861-6878.

\bibitem{Fer4} E.V.~Ferapontov, Dupin hypersurfaces and integrable
Hamiltonian systems of hydrodynamic type which do not possess Riemann
invariants, Diff. Geom. and its Appl. 5 (1995) 121-152.


\bibitem{Finikov} S.P.~Finikov, Theory of congruences, M-L.,
Gostekhizdat, 1950.


\bibitem{Forsyth} A.R.~Forsyth, Theory of differential equations,
Cambridge University Press, 6 (1906).


\bibitem{Kamran} N.~Kamran and K.~Tenenblat,
Laplace transformations in higher dimensions, Duke Math. J. 1996. V.84, N1
P.237-266.


\bibitem{Levy} L\'evy, Journ. de l' \'Ec. Polytechnique, t. XXXVII,
Cah. LVI (1886).

\bibitem{Sevennec} B.~S\'evennec, G\'eom\'etrie des syst\`emes
hyperboliques
de lois de conservation, M\'emoire (nouvelle s\'erie) N56,
Suppl\'ement au Bulletin de la Soci\'et\'e Math\'ematique de France,
V. 122, 1994, 1-125.

\bibitem{Tsarev} S.P.~Tsarev, The geometry of Hamiltonian systems of
hydrodynamic type. The generalized hodograph transform, Math. USSR Izv.,
37 (1991) 397-419.


\end{thebibliography}
\end{document}